\documentclass{amsart}
\usepackage{amsmath,amssymb}
\usepackage{bbm}
\usepackage{graphicx,tikz,mathtools}

\usepackage{color}
\usepackage{comment}
\definecolor{Caca}{RGB}{193,124,250}

\newtheorem{theorem}{Theorem}[section]

\theoremstyle{definition}

\theoremstyle{remark}

\newcommand{\Si}{\Sigma}

\newcommand\minus\backslash

\newcommand\lan\langle
\newcommand\ran\rangle

\renewcommand\leq\leqslant
\renewcommand\geq\geqslant
\newlength{\intwidth}

\numberwithin{equation}{section}

\newcommand\BSi\Si 

\begin{document}

\title[]{Absence of uniform bounds for  the boundary nodal count of Steklov eigenfunctions on Euclidean domains}

 \author{Alberto Enciso}
 \address{Instituto de Ciencias Matem\'aticas, Consejo Superior de
   Investigaciones Cient\'\i ficas, 28049 Madrid, Spain, e-mail: {\sf aenciso@icmat.es}.}

\author{Angela Pistoia}
\address{Dipartimento di Scienze di Base e Applicate per l'Ingegneria, Universit\`a di Roma ``La Sapienza'',  00161 Roma, Italy, e-mail: {\sf angela.pistoia@uniroma1.it}.}

\author{Luigi Provenzano}
\address{Dipartimento di Matematica e Informatica, Universit\`a degli Studi di Cagliari,  09124 Cagliari, Italy, e-mail: {\sf luigi.provenzano@unica.it}.}

%
\begin{abstract}
 We prove that for any $K,N\in\mathbb N$ there exist smooth bounded domains in $\mathbb R^n$, $n\geq 3$, homeomorphic to a ball, for which the first $N$ positive Steklov eigenvalues are simple and the boundary traces of the corresponding eigenfunctions have at least $K$ nodal domains. This provides a negative answer to the question of whether the number of nodal domains of Dirichlet-to-Neumann eigenfunctions satisfies a Courant-type bound in the Euclidean setting and complements the results of \cite{EPP}.
\end{abstract}

\maketitle

\section{Introduction and statement of the main result}

Let $\Omega$ be a smooth  bounded domain in $\mathbb R^n$. The Dirichlet-to-Neumann map $\mathcal D:H^{1/2}(\partial\Omega)\to H^{-1/2}(\partial\Omega)$ is defined as $\mathcal D\phi=\partial_\nu u$, where $u$ is the {\it harmonic extension} of $\phi$, namely
$$
\Delta u=0\ {\rm in\ }\Omega\,,\ \ \ u=\phi\ {\rm on\ }\partial\Omega.
$$
 It is well known that $\mathcal D$ is a nonlocal pseudodifferential elliptic operator on $\partial\Omega$
of order 1 which defines a non-negative, self-adjoint operator with dense domain in $L^2(\partial\Omega)$. In particular, one can find an orthonormal basis of $L^2(\partial\Omega)$ of Dirichlet-to-Neumann eigenfunctions $\{\phi_k\}_{k=0}^\infty$ satisfying the eigenvalue equation \begin{equation}\label{DtN}
\mathcal D\phi_k=\sigma_k\phi_k.
\end{equation}
The harmonic extensions of $\phi_k$ to $\Omega$, call them $u_k$, are Steklov eigenfunctions, as they solve the Steklov eigenvalue problem
\begin{equation}\label{Steklov}
\begin{cases}
\Delta u_k=0\,, & {\rm in\ }\Omega,\\
\partial_{\nu}u_k=\sigma_ku_k & {\rm on\ }\partial\Omega.
\end{cases}
\end{equation}
If $\Omega$ is merely Lipschitz, problems \eqref{DtN} and \eqref{Steklov} are understood in the weak sense. The eigenvalues are non-negative, of finite multiplicity, and form an increasing sequence diverging to $+\infty$
$$
0=\sigma_0(\Omega)<\sigma_1(\Omega)\leq\cdots\leq\sigma_k(\Omega)\leq\cdots\nearrow+\infty.
$$
We  call them Steklov or Dirichlet-to-Neumann eigenvalues and highlight their dependence on the domain $\Omega$. Note that $u_0$ is constant, while $u_1$ changes sign in $\Omega$, and accordingly, $\phi_0$ is constant on $\partial\Omega$, $\phi_1$ changes sign on $\partial\Omega$. We refer to the  surveys by Girouard and Polterovich~\cite{Pol17} and by Colbois, Girouard, Gordon and Sher~\cite{GirouSurv} for an exhaustive treatment of the Steklov and Dirichlet-to-Neumann spectral problems.

In this short note we are interested in the {\it nodal sets} of Dirichlet-to-Neumann eigenfunctions $\phi_k$, that is, the nodal sets of the restriction of Steklov eigenfunctions $u_k$ to the boundary.

Courant's nodal domain theorem bounds the number of nodal domains of $u_k$ in terms of $k$: they are at most $k+1$, see \cite{hash,KS}. By boundary nodal domains we intend the nodal domains of the Dirichlet-to-Neumann eigenfunctions $\phi_k=u_k|_{\partial\Omega}$, that is, the connected components of $\partial\Omega\setminus\phi_k^{-1}(0)$. In dimension $2$, bounds  on interior and boundary nodal domains are related by the topology, and in fact the number of nodal domains of $\phi_k$, $k\geq 1$, is at most $2(k+b)$  for any planar domain with $b+1$ boundary components. A  question that has been featured in the surveys \cite[Open problem 9]{Pol17} and \cite[Open question 10.14]{GirouSurv} is whether for $n\geq 3$ the number of boundary nodal domains satisfies a Courant-type upper bound, possibly depending on the topology.

 In \cite{EPP} we showed that, given a compact Riemannian manifold $M$ with boundary and dimension $n\geq 3$, and given any $K,N\in\mathbb N$ (where by $\mathbb N$ we intend the set of {\it positive} natural numbers), it is always possible to find a metric  $g$ for which $\phi_k$ has at least $K$ nodal domains, for all $k=1,...,N$.  This happens regardless of the topology of $M$ which can be homeomorphic to a ball. We left open the question whether there exist {\it Euclidean domains}, in dimension at least $3$, for which the same conclusion holds, that is, for which $\phi_k$ has at least $K$ nodal domains for all $k=1,...,N$. In this note we address this question and provide  a simple construction of Euclidean domains, homeomorphic to balls, for which $\sigma_1,...,\sigma_N$ are simple and $\phi_1,...,\phi_N$ have at least $K$ nodal domains, for any choice of $K,N\in\mathbb N$. We remark that in \cite{EPP} we had the freedom of moving a metric on a fixed manifold. In the case of domains the situation is more rigid and the same approach cannot work. Instead, in our Euclidean examples, the complexity of the boundary nodal domains comes entirely from the shape of the domain.

 Here is the  result:

\begin{theorem}\label{mainT}
Let $n\geq 3$. For any $K,N\in\mathbb N$ there exists a smooth bounded domain in $\mathbb R^n$  homeomorphic to a ball such that the first $N$ positive eigenvalues are simple and the corresponding Dirichlet-to-Neumann eigenfunctions have at least $K$ nodal domains.
\end{theorem}

{\bf Sketch of the construction for $N=1$.} In order to understand the main ideas behind the proof of Theorem \ref{mainT}, we find it convenient to sketch it for the first non-trivial eigenfunction, that is, for $N=1$. Consider the union of the balls $B(0,1)$, $B(p_1(\epsilon),R)$ and $B(p_2(\epsilon),R)$, where $p_1(\epsilon)=(-1-R+\epsilon,0,...,0)$ and $p_2(\epsilon)=(1+R-\epsilon,0,...,0)$. This union is connected and Lipschitz for $\epsilon>0$ small. Then, remove $K$ spherical holes with centers in $\{x_1=0\}$ from $B(0,1)$. Here $x=(x_1,...,x_n)$ denotes a point in $\mathbb R^n$. The resulting domain has a boundary made of $K+1$ connected components.  See Figure \ref{fig1} for a graphic representation.
As $\epsilon\to 0$ the domain eventually degenerates to the disjoint union of the ``central'' ball $B(0,1)$ with  $K$ spherical holes, and $2$ tangent ``lateral'' balls $B(p_1(0),R)$, $B(p_2(0),R)$. The whole domain is invariant under the reflection $x_1\mapsto -x_1$. The symmetry of the domain forces the first non-trivial Steklov eigenspace to be spanned by even and odd eigenfunctions in $x_1$. However, as $\epsilon\to 0$, the eigenfunctions converge to constants in each ball. At this point, when $\epsilon$ is small enough,  Courant's nodal domain (which holds for Steklov eigenfunctions, i.e., in the interior) prevents a first non-trivial eigenfunction to be even in $x_1$, and also prevents the existence of two linearly independent odd eigenfunctions. Hence a first Steklov eigenvalue must be simple. The odd symmetry and Courant's theorem force the eigenfunction to have exactly two interior nodal domains, contained respectively in $\{x_1>0\}$ and $\{x_1<0\}$. This implies  two boundary nodal domains on each of the $K$ holes and two on the outer boundary. Overall, $2K+2$ boundary nodal domains. Well-known perturbation results imply that we can approximate this domain with smooth domains. Importantly, the fact that there are at least two nodal domains on the boundary of each (possibly perturbed) hole has not changed. To conclude, it remains to change the topology of the domain without losing too many boundary nodal domains. We do so through {\it tubular excisions}: we connect the holes to the outer boundary by removing thin tubes around smooth embedded curves: this produces a domain homeomorphic to a ball, and if the size of the tubes is small enough, we have a nice convergence of the eigenfunctions up to the boundary far from the curves. This keeps at least $K+1$ of the boundary nodal domains of the original set. See Figure \ref{fig2} for a graphic representation of the domain after the surgery.

{\bf The general case.} For higher eigenfunctions, that is, when $N\geq 1$, we reproduce a similar situation: a central ball with $K$ holes with $N+1$ symmetrically placed lateral balls. More precisely, we consider connected, Lipschitz domains, invariant for the action of the cyclic group of rotations of integer multiples of $\frac{2\pi}{N+1}$ on a fixed plane $\Sigma$, from which $K$ spherical holes with centers in $\Sigma^\perp$ are removed. The domain eventually degenerates to the disjoint union of a ``central'' ball with  $K$ spherical holes, and $N+1$ ``lateral'' balls placed symmetrically around the central ball (their centers lie on a circle in $\Sigma$ and are equally spaced).  Figure \ref{fig1} is indeed a representation of the case $N=1, K=3$. The symmetry of the domain and the placement of the holes force the first $N$ non-constant eigenfunctions to change sign on the boundary of each of the $K$ holes, which means that there are at least $2K+1$ boundary nodal domains. At this point, the remaining arguments are the same as in the case $N=1$. To simplify the application of the perturbation results, we want the first $N$ eigenvalues to be simple, which we may always assume by generic simplicity.

The use of tubular excisions for the  Steklov problem in higher dimensions is not new, and has been used, for example, to produce examples of large eigenvalues \cite{bri1,bri2}, examples of domains homeomorphic to a ball having larger normalised first eigenvalue than that of the ball \cite{FS_SO}, and in general to study spectral optimisation problems for higher eigenvalues \cite{hong}. The idea is that of transferring information which is easier to get on domains with non-trivial topologies (holes, homogenisation, etc.) to domains with simpler topology.

The proof of Theorem \ref{mainT} is contained in Section \ref{sec:proof} here below. We have also included Subsection \ref{sub:sigma1} where, for the reader's convenience, we have highlighted the simpler proof for the first eigenfunction sketched here above.  In Appendix \ref{app:A} we have collected a few known results on singular perturbation and spectral stability for the Steklov problem that we have employed in the proof of Theorem \ref{mainT}.

\section{Proof of Theorem \ref{mainT}}\label{sec:proof}
{\bf Step 1. A symmetric Lipschitz domain with spherical holes.}
Let $(e_1,...,e_n)$ be an orthonormal basis of $\mathbb R^n$, and let $(x_1,...,x_n)$ be the corresponding Cartesian coordinates. We write $x=(x_1,...,x_n)$ for $x\in\mathbb R^n$. By $B(p,r)$ we denote a ball in $\mathbb R^n$ centered at $p$ and of radius $r$.

\begin{figure}
    \centering
    \includegraphics[width=0.8\textwidth]{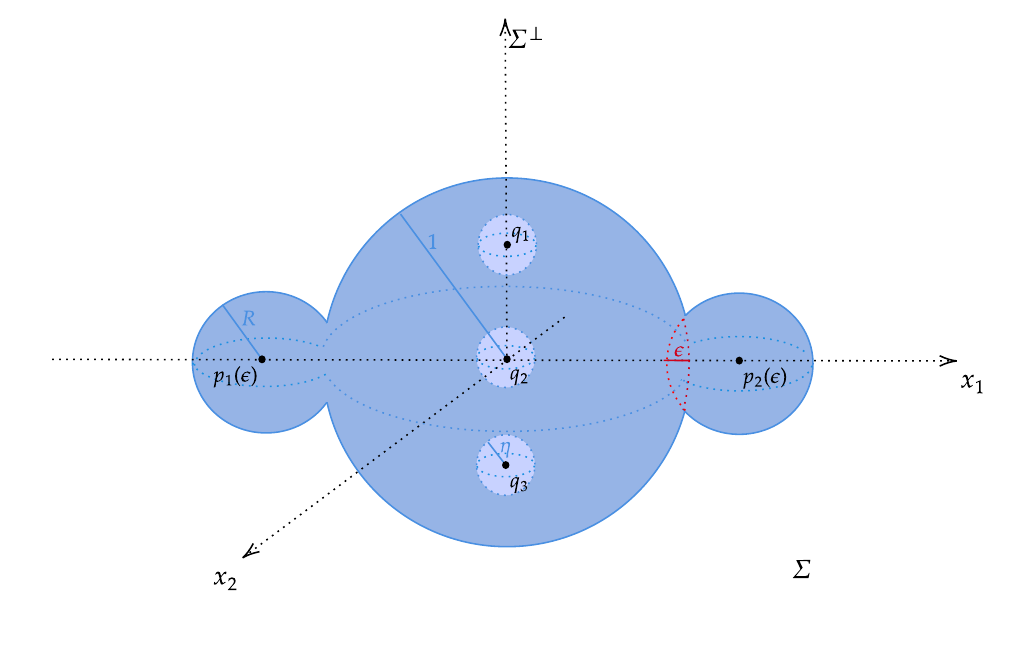}
    \caption{The domain $\Omega_\epsilon$ with $N=1$ and $K=3$.}
    \label{fig1}
\end{figure}

Let $\Sigma:=\{x\in\mathbb R^n:x_3=x_4=\cdots=x_n=0\}$ be the plane in $\mathbb R^n$ spanned by $e_1,e_2$, and let $\Sigma^\perp$ be its orthogonal complement, that is $\Sigma^\perp:=\{x\in\mathbb R^n:x_1=x_2=0\}$.  Fix $N\in\mathbb N$. For any $\epsilon\in[0,1]$ consider then the following $N+1$ points in $\Sigma$:
$$
p_i(\epsilon)=(1+R-\epsilon)(\cos(\frac{2\pi i}{N+1}),\sin(\frac{2\pi i}{N+1}),0,...,0)\,,\ \ \ i=1,...,N+1.
$$
Fix $R>0$ such that the balls $B(p_i(\epsilon),R)$ are pairwise disjoint for all $\epsilon\in[0,\epsilon_0)$, for some $\epsilon_0>0$. Let $K\in\mathbb N$ and consider $K$ points on $\Sigma^\perp$, $q_1,...,q_K$ and $\eta>0$ such that the $K$ closed balls $\overline{B(q_i,\eta)}$  are pairwise disjoint and $\overline B(q_i,\eta)\subset B(0,1)$. All the parameters $N,K,R,\eta$ are fixed, and we let just $\epsilon$ vary. Set, for $\epsilon\in[0,\epsilon_0)$:
$$
\Omega_\epsilon:=\left(B(0,1)\bigcup_{i=1}^{N+1}B(p_i(\epsilon),R)\right)\setminus\bigcup_{i=1}^K \overline {B(q_i,\eta)}
$$
See Figure \ref{fig1} for a representation in the case $N=1,K=3$.

By construction, possibly choosing a smaller value of $\epsilon_0$, we have that for all $\epsilon\in(0,\epsilon_0)$, $\partial B(p_i(\epsilon),R)\cap \partial B(q_i,\eta)= \emptyset$ for all $i=1,...,N+1$, and $\Omega_\epsilon$ is  piecewise smooth and Lipschitz. It degenerates to $\Omega_0$, which is the  disjoint union of the $N+1$ balls $B(p_i(0),R)$ and of the domain with holes $\Omega_c:=B(0,1)\setminus \bigcup_{i=1}^K \overline B(q_i,\eta)$.

The well-known spectral stability under degeneration of a domain to a disjoint union of domains, which we have stated in Theorem \ref{thm:degeneration}, implies that 
$$
\sigma_k(\Omega_\epsilon)\to 0
$$
for all $k=1,...,N+1$ and $\sigma_{N+2}\to\min\{\sigma_1(\Omega_c),\frac{1}{R}\}$. See \cite{bucur_giacomini_trebeschi,hong} for the relevant references.  In fact, Theorem \ref{thm:degeneration} states that the spectrum of 
$\Omega_\epsilon$ converges to the spectrum of
$$\Omega_0:=\Omega_c\sqcup B(p_1(0),R)\sqcup \cdots\sqcup B(p_{N+1}(0),R),
$$
where by $\sqcup$ we denote the disjoint union of sets.
The spectrum of this disjoint union is just the union of the spectra of $N+1$ copies of the ball of radius $R$ and of $\Omega_c$.

{\bf Step 2. The first $N$ Steklov eigenfunctions change sign  on the boundary of each hole. }We prove now that any Dirichlet-to-Neumann eigenfunction $\phi_\epsilon$ associated with $\sigma_k(\Omega_\epsilon)$ has at least $2K+1$ nodal domains for all $k=1,...,N$. Denote by $R_{N+1}$ the rotation of $\frac{2\pi}{N+1}$ in the plane $\Sigma$ about the origin. By construction, $\Omega_\epsilon$ is invariant under the action of $R_{N+1}$. 

First, we show that for $\epsilon$ small no such $\phi_\epsilon$ can satisfy $\phi_\epsilon\circ R_{N+1}=\phi_\epsilon$. If this were not true, by Theorem \ref{thm:degeneration}, we would find a sequence $\epsilon_i\to 0$ as $i\to\infty$ such that $\phi_{\epsilon_i}\to\phi_0$ in $C^0(K)$  as $i\to\infty$ for all compact $K$ of $\overline\Omega_0\setminus\{w_1,...,w_{N+1}\}$, where $w_i=R_{N+1}^ie_1$ is the point of tangency of $\overline{B(p_i(0),R)}$ with $\overline\Omega_c$ and $\phi_0$ is some eigenfunction associated with $\sigma_k(\Omega_0)$. Here we write $R_{N+1}^i:=R_{N+1}\circ\cdots\circ R_{N+1}$ for the $i$-th power of $R_{N+1}$. Note that, since $k\leq N+1$, $\phi_0$ is constant on each connected component of $\Omega_0$ (the constant can be zero on some connected components). The functions $\phi_\epsilon$, $\phi_0$ are normalised by $\int_{\partial\Omega_{\epsilon_i}}\phi_{\epsilon_i}^2=1$, $\int_{\partial\Omega_0}\phi_0^2=1$. Then, also  $\phi_0$ satisfies $\phi_0\circ R_{N+1}=\phi_0$, which means that $\phi_0=a$ on $B(p_i(0),R)$, $i=1,...,N+1$ and $\phi_0=b$ on $\Omega_c$, for some $a,b\in\mathbb R$ which must satisfy  $\int_{\partial\Omega_0}\phi_0=0$, that is $(N+1)|\partial B(0,R)|a+|\partial\Omega_c|b=0$  and from the normalisation, $(N+1)|\partial B(0,R)|a^2+|\partial\Omega_c|b^2=1$ which forces $a,b\ne 0$, $ab<0$. This implies that for $\epsilon_i$ small enough $\phi_{\epsilon_i}$ has at least $N+2$ nodal domains on the outer connected component of the boundary. Consider then the Steklov eigenfunction $u_{\epsilon_i}$ (harmonic extension of $\phi_{\epsilon_i}$).

By the $C^0$ convergence up to the boundary (up to the choice of signs) $u_{\epsilon_i}$ is negative on the closure of the central component outside sufficiently small, fixed and pairwise disjoint neighbourhoods of the intersections $\Gamma_i:=\partial\Omega_c\cap\partial{B(p_i(\epsilon),R)}$, and positive at a point in each lateral ball $B(p_i(\epsilon),R)$, for $i$ sufficiently large. Since any path joining two lateral balls must cross the region where $u_{\epsilon_i}$ is negative, these positive points belong to distinct interior nodal domains. Together with a negative nodal domain in the central component, this gives at least $N+2$ interior nodal domains, contradicting Courant's nodal domain theorem for Steklov eigenfunctions \cite{hash,KS} (i.e., in the interior) for $k\leq N$.

Since $N$ is fixed, we find $\epsilon>0$ so that no eigenfunction $\phi_\epsilon$  associated with some $\sigma_1(\Omega_\epsilon),..,\sigma_N(\Omega_\epsilon)$ can be invariant by the action of $R_{N+1}$, that is, there is no eigenfunction $\phi_\epsilon$ such that $\phi_\epsilon\circ R_{N+1}=\phi_\epsilon$. From now on we fix this $\epsilon$.

Let $\phi$ be some eigenfunction associated with $\sigma_k(\Omega_\epsilon)$. Note that $\phi\circ R_{N+1}^i$ is still an eigenfunction in the same eigenspace. Take
$$
A\phi:=\frac{1}{N+1}\sum_{i=1}^{N+1}\phi\circ R_{N+1}^i
$$
Then $A\phi\circ R_{N+1}=A\phi$ and by the discussion above, $A\phi=0$.
Now, all the removed holes $B(q_i,\eta)$ are invariant under $R_{N+1}$, that is, $R_{N+1}(B(q_i,\eta))=B(q_i,\eta)$, hence
$$
\int_{\partial B(q_i,\eta)}\phi=\int_{\partial B(q_i,\eta)}A\phi=0.
$$
Since $\phi$ is a non-constant Steklov eigenfunction, by unique continuation we have that $\phi|_{\partial B(q_i,\eta)}$ cannot be identically zero, hence it has to change sign. Therefore $\phi$ has at least $2K+1$ nodal domains, of which at least two on each $\partial B(q_i,\eta)$.

{\bf Step 3. Smoothing out the set and making the first $N$ eigenvalues simple.} The domain $\Omega_\epsilon$ is Lipschitz and piecewise smooth, and it is invariant under the action of $R_{N+1}$. Hence each eigenvalue $\sigma_1(\Omega_\epsilon),...,\sigma_N(\Omega_\epsilon)$ may have multiplicity.

Now, we approximate $\Omega_\epsilon$ with smooth domains in the following way. Since $\Omega_\epsilon$ is smooth far from the intersections $\Gamma_i=\partial\Omega_c\cap\partial{B(p_i(\epsilon),R)}$, it is enough to smooth out the domain in a neighborhood of each $\Gamma_i$. Note that in a neighborhood of $\Gamma_i$,  $\Omega_\epsilon$ is just the union of two euclidean balls. For simplicity, we describe the smoothing of the union of two balls $B_1,B_2\subset\mathbb R^n$ with non-empty intersection, which is the kind of regularisation that we then apply to $\Omega_\epsilon$. After a rescaling and a congruence, we may suppose that
$$
\partial B_1\cap \partial B_2=\{(x_1,...,x_n)\in\mathbb R^n:x_n=0, x_1^2+\cdots+x_{n-1}^2=1\}.
$$
In this coordinate system
$$
B_1=\{(x_1,...,x_n):x_1^2+\cdots+x_n^2<1+2ax_n\}$$
$$
B_2=\{(x_1,...,x_n):x_1^2+\cdots+x_n^2<1+2bx_n\}
$$
for some $a<b$. Then
$$
B_1\cup B_2=\{(x_1,...,x_n):x_1^2+\cdots+x_n^2<1+(a+b)x_n+(b-a)|x_n|\}.
$$
Define, for $\delta\in(0,\delta_0)$, 
$$
D_\delta:=\{(x_1,..,x_n):x_1^2+\cdots+x_n^2<1+(a+b)x_n+(b-a)\sqrt{x_n^2+\delta^2\chi(x_n/\delta)}\}
$$
where $\chi\in C^{\infty}_c((-1,1))$, $0\leq\chi\leq 1$  and $\chi\equiv 1$ in a neighborhood of $0$. We check that $D_\delta\to B_1\cup B_2$ (for the Hausdorff distance) and $|\partial D_\delta|\to |\partial(B_1\cup B_2)|$ as $\delta\to 0$ (convergence of the boundary measures). Setting $g_0(x_n):=1+(a+b)x_n+(b-a)|x_n|$ and $g_\delta(x_n):=1+(a+b)x_n+(b-a)\sqrt{x_n^2+\delta^2\chi(x_n/\delta)}$ (the functions defining the radial profiles of $B_1\cup B_2$ and $D_\delta$), we also see that $g_\delta\to g_0$ in $L^\infty(\mathbb R)$ and $\|g_\delta'\|_{L^\infty(\mathbb R)}<C$ independently of $\delta$. We are exactly in the classical hypotheses of the spectral stability for the Steklov problem \cite{FerreroLamberti}, which we have summarised in Theorem \ref{thm:spec:stab}. The assumptions of Theorem \ref{thm:spec:stab} are satisfied by the family $\{D_\delta\}_{\delta\in(0,\delta_0)}$, hence we have convergence of the spectrum to that of $B_1\cup B_2$ as $\delta \to 0$ and, up to subsequences, of the eigenfunctions in $C^0$ away from a small neighborhood of $\partial B_1\cap\partial B_2$. We apply the smoothing procedure described above to $\Omega_\epsilon$ near each intersection $\Gamma_i$. Doing so, we get a smooth domain $\Omega_{\epsilon,\delta}$ which coincides with $\Omega_\epsilon$ outside a small neighborhood of each $\Gamma_i$. The domain $\Omega_{\epsilon,\delta}$ is smooth, and by Theorem \ref{thm:spec:stab} we have that any eigenfunction associated with $\sigma_k(\Omega_{\epsilon,\delta})$ has at least $2K+1$ nodal domains, of which at least two on each $\partial B(q_i,\eta)$. In fact, if this were not the case, Theorem \ref{thm:spec:stab} would provide a sequence of normalized eigenfunctions failing to change sign on a fixed hole boundary and converging, up to a subsequence, uniformly on the boundary of the hole to a normalized eigenfunction of $\Omega_\epsilon$ with the same property, contradicting Step 2.

Note that at this stage, the eigenvalues may be still multiple. Here we invoke a generic simplicity result for Steklov eigenvalues \cite{wang_generic_euclid}, which we  have recalled in Theorem \ref{thm:genericity}.
From Theorem \ref{thm:genericity} we deduce that there exists a smooth embedding $\Phi:\overline\Omega_{\epsilon,\delta}\to\mathbb R^n$ such that $\sigma_j(\Phi(\Omega_{\epsilon,\delta}))$, $j=1,...,N$ are simple. This map can be chosen arbitrarily close to the identity in the $C^m(\Omega_{\epsilon,\delta};\mathbb R^n)$ topology, for all $m\geq 2$, and to be the identity on $\overline\Omega_{\epsilon,\delta}\cap U$, where $U$ is any prescribed open subset of $\mathbb R^n$ intersecting $\partial\Omega_{\epsilon,\delta}$ and such that $\partial\Omega_{\epsilon,\delta}\setminus\overline U\ne\emptyset$. We choose $\Phi$ to be the identity on a neighborhood of the holes $B(q_i,\eta)$ and to be close enough to the identity in $C^m(\overline\Omega_{\epsilon,\delta};\mathbb R^n)$ for some $m\geq 2$ so that the spectral stability theorem (Theorem \ref{thm:spec:stab}) applies. Then, the resulting domain $\Phi(\Omega_{\epsilon,\delta})$ is of class $C^m$, $m\geq 2$, the first $N$ positive eigenvalues are simple and the corresponding eigenfunctions have at least $2K+1$ nodal domains, of which at least two on each $\partial B(q_i,\eta)$ (which have not been touched by any of the perturbation above). We set now $\tilde\Omega:=\Phi(\Omega_{\epsilon,\delta})$.

\begin{figure}
     \centering
     \includegraphics[width=0.8\textwidth]{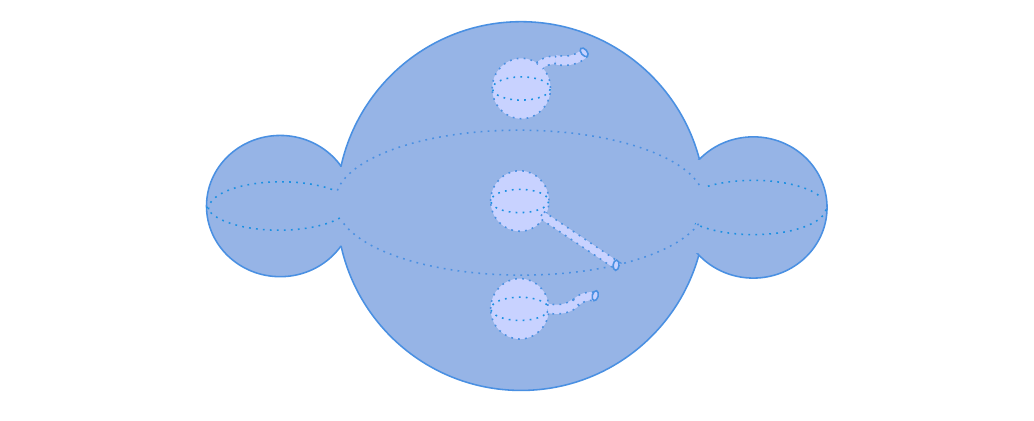}
     \caption{The domain $\Omega_\epsilon$ after tubular excisions.}
     \label{fig2}
 \end{figure}

{\bf Step 4. } By construction $\tilde\Omega$ is homeomorphic to a ball with $K$ disjoint balls removed from its interior. The next step is to perturb $\tilde \Omega$ in order to make it homeomorphic to a ball. Let $\sigma_1(\tilde\Omega),...,\sigma_N(\tilde\Omega)$ be the first $N$ positive Steklov eigenvalues and call $u_k$ the corresponding eigenfunctions, normalised by $\int_{\partial\tilde\Omega}u_k^2=1$. They are unique up to the choice of the sign. Choose $K$ points $y_i\in\partial B(q_i,\eta)$ such that  $u_k(y)\ne 0$ for all $y\in\overline D_i\subset\partial B(q_i,\eta)$, for all $k=1,...,N$  , where $D_i$ is some open geodesic disk centered at $y_i$ on $\partial B(q_i,\eta)$. This is possible by unique continuation. Then, for each $k=1,...,N$ choose points $\xi_{i,k}\in\partial B(q_i,\eta)$ such that $u_k(\xi_{i,k})u_k(y_i)<0$. Let also $A_i=\partial D_i\subset\partial B(q_i,\eta)$. By construction $u_k(y)\ne 0$ on $A_i$ and $u_k(y)u_k(\xi_{i,k})<0$ for all $y\in A_i$.
 Choose distinct points $z_1,...,z_K\in\partial\tilde\Omega\setminus\bigcup_{i=1}^K\partial B(q_i,\eta)$ (i.e., points on the outer boundary) and smooth, disjoint, embedded curves $\gamma_i:[0,1]\to\mathbb R^n$, $\gamma_i(0)=y_i$, $\gamma_i(1)=z_i$, $\gamma_i((0,1))\subset\tilde\Omega$, meeting $\partial\tilde\Omega$ orthogonally at the extrema. Let $\Gamma=\bigcup_{i=1}^K\gamma_i$ and define, for $\rho\in(0,\rho_0)$ sufficiently small,
 $$
\tilde\Omega_\rho:=\{x\in\tilde\Omega:{\rm dist}(x,\Gamma)>\rho\}.
 $$
 See figure \ref{fig2} for a graphic representation.
 
 This domain is homeomorphic to a ball: it is obtained from $\tilde\Omega$ by removing thin tubes around smooth embedded curves connecting the holes with the outer boundary. In fact, each excision joins two distinct spherical boundary components by a connected sum, so that the resulting boundary is again a sphere. It is well-known that spectral stability  holds also under this singular perturbation, when $n\geq 3$, see \cite{FS_SO,hong}. We have stated the precise result in Theorem \ref{thm:excisions}. From Theorem \ref{thm:excisions} we have that $\sigma_k(\tilde\Omega_\rho)$, $k=1,...,N$ are simple for all $\rho\in(0,\rho_0)$ (possibly choosing a smaller $\rho_0$). Moreover, if $u$ denotes an eigenfunction associated with $\sigma_k(\tilde\Omega)$ normalised by $\int_{\partial\tilde\Omega}u^2=1$ and $u_\rho$ denotes an eigenfunction associated with $\sigma_k(\tilde\Omega_\rho)$ normalised by $\int_{\partial\tilde\Omega_\rho}u_\rho^2=1$, Theorem \ref{thm:excisions} implies that (up to the choice of signs) $u_\rho\to u$ in $C^0(K)$ for all compact $K$ in the closure of $\tilde\Omega$ away from the set $\Gamma$ of the limit curves. 
 
 In particular, $A_i\subset\partial\tilde\Omega_\rho\cap \partial B(q_i,\eta)$ and $u_\rho(y)u_\rho(\xi_{i,k})<0$ for all $y\in A_i$ when $\rho>0$ is sufficiently small. The boundary nodal domain containing $\xi_{i,k}$ cannot cross $A_i$, where $u_\rho$ has the opposite sign, and is therefore contained in $\partial B(q_i,\eta)\setminus\overline{D_i}$. These give $K$ distinct boundary nodal domains. Since a point of some $A_i$ belongs to an additional boundary nodal domain, overall there are at least $K+1$ boundary nodal domains. The domain $\tilde\Omega_\rho$ is piecewise smooth and Lipschitz. It has corners where the tubes around $\gamma_i$ meet $\partial\tilde\Omega_\rho$. Smoothing the corners as in Step 3, while leaving neighborhoods of the sets $A_i$ and the points $\xi_{i,k}$ unchanged, yields a smooth domain $\Omega$, homeomorphic to a ball. By Theorem \ref{thm:spec:stab}, the first $N$ positive eigenvalues remain simple and each corresponding eigenfunction has at least $K+1$ boundary nodal domains. This concludes the proof. \qed

\subsection{Remark: a simple argument for $u_1$}\label{sub:sigma1} The proof assumes a simpler form if we confine to the basic situation of $u_1$. We sketch it here since we think it can be useful. Consider, for all $\epsilon\in(0,1/2)$, the domain
  $$
\Omega_\epsilon:=B(-p(\epsilon),1)\cup B_K\cup B(p(\epsilon),1)
  $$
  where
  $$
p(\epsilon)=(2-\epsilon,0,...,0)
  $$
   and
  $$
B_K:=B(0,1)\setminus \bigcup_{i=1}^K \overline{B(q_i,\eta)}
  $$
where $q_i\in\{(x_1,...,x_n):x_1=0\}$ and $\overline{B(q_i,\eta)}$ are pairwise disjoint and strictly contained in $B(0,1)$. In other words: $B_K$ is the unit ball with $K$ holes with centers on $\{x_1=0\}$. By  symmetry, the first eigenspace admits a basis of eigenfunctions that are either even or odd in $x_1$. However, as $\epsilon\to 0$ a first eigenfunction cannot be even. In fact, thanks to Theorem \ref{thm:degeneration}, up to subsequences, it would converge (up to signs) to a constant $a>0$ on $B(-p(0),1)\cup B(p(0),1)$ and a constant $b<0$ on $B_K$. In particular, it would converge to  $b$ in $\overline B_K\cap\{-1/2<x_1<1/2\}$. This would imply at least three interior nodal domains for a subsequence of first eigenfunctions, violating Courant's nodal domain theorem in the interior. Hence, for all $\epsilon>0$ small, a first eigenfunction $u_1$ must be odd in $x_1$ and the corresponding eigenspace of dimension $1$, again by Courant's nodal domain theorem:  two odd, orthogonal first eigenfunctions $u,u'$ would satisfy $uu'>0$ or $uu'<0$ on $\overline\Omega_\epsilon\setminus\{x_1=0\}$ and also $\int_{\partial\Omega_\epsilon}uu'=0$, which is impossible. By symmetry, for all $\epsilon>0$ small (which we may fix from now on), $u_1$ changes sign on the boundary of any hole $B(q_i,\eta)$. Smoothing and tubular excisions as in Step 4 of the proof provide a domain homeomorphic to a ball with at least $K+1$ boundary nodal domains (which can be again approximated by smooth domains, resulting in a smooth example).

  \appendix

  \section{Known results on singular perturbation and spectral stability}\label{app:A}

In this appendix we collect some classical results in spectral perturbation for the Steklov problem. 

 The first result is a convergence result for Lipschitz domains degenerating to the disjoint union of smooth domains. Just to fix the ideas, it is well-known that for the union of two disks  in the plane (with non-empty intersection) degenerating to the disjoint union of two tangent disks, the Steklov spectrum converges to the union of the spectra of the two disks (see e.g., \cite{Pol17} and references therein). The same phenomenon happens in more generality: it holds in any dimension and when the domain degenerates to the disjoint union of $N$ smooth domains tangent at a finite number of elliptic points. This is done in \cite[Theorem 3.2 and Remark 3.3]{hong}, where also the convergence of eigenfunction is discussed. We state the result in a specific situation of a  ``central'' domain and $N$ attached ``lateral'' domains.
\begin{theorem}\label{thm:degeneration}
Let $\Omega_0,...,\Omega_N$ be bounded smooth domains in $\mathbb R^n$, $n\geq 2$, such that $\Omega_0\cap\Omega_j=\emptyset$, $\overline\Omega_0\cap\overline\Omega_j=\{p_j\}$ for all $j=1,...,N$ and $\overline\Omega_i\cap\overline\Omega_j=\emptyset$ for all $i\ne j=1,...,N$ and the principal curvatures of $\partial\Omega_0$ and $\partial\Omega_j$ at $p_j$ are all positive (with the convention that the boundary of a ball has positive principal curvatures). Let $e_j\in\mathbb S^{n-1}$ and $\epsilon_0>0$ be such that $(\epsilon e_j+\Omega_j)\cap\Omega_0\ne\emptyset$ for all $j=1,...,N$ and $\overline{(\epsilon e_j+\Omega_j)}\cap\overline{(\epsilon e_i+\Omega_i)}=\emptyset$ for all $i\ne j=1,...,N$ and $\epsilon\in(0,\epsilon_0)$ and $\Omega_\epsilon:=\Omega_0\bigcup_{j=1}^N(\epsilon e_j+\Omega_j)$ is connected for all $\epsilon\in(0,\epsilon_0)$. Then
\begin{equation}\label{spectral_stab}
\lim_{\epsilon\to 0}\sigma_k(\Omega_\epsilon)=\sigma_k(\Omega_0\sqcup\Omega_1\sqcup\cdots\sqcup\Omega_N),
\end{equation}
where $\sqcup$ denotes the disjoint union of sets. Moreover, let $\{u_{\epsilon}\}_{\epsilon\in(0,\epsilon_0)}$ be a family of Steklov eigenfunctions corresponding to $\sigma_k(\Omega_\epsilon)$, normalised by $\int_{\partial\Omega_\epsilon}u_{\epsilon}^2=1$. Then, for any sequence $\epsilon_i\to 0$ as $i\to\infty$ there exists a subsequence $\epsilon_{i_\ell}\to 0$ as $\ell\to\infty$ and a Steklov eigenfunction $u$ associated with $\sigma_k(\Omega_0\sqcup\cdots\sqcup\Omega_N)$ normalised by $\int_{\partial\Omega_0\sqcup\cdots\sqcup\partial\Omega_N}u^2=1$ such that
\begin{equation}\label{spectral_stab_eig_1}
\lim_{\ell\to\infty}\|u_{\epsilon_{i_\ell}}(\cdot+\epsilon_{i_\ell} e_j)-u\|_{H^1(\Omega_j)}=0
\end{equation}
and
\begin{equation}\label{spectral_stab_eig_2}
\lim_{\ell\to \infty}\|u_{\epsilon_{i_\ell}}(\cdot+\epsilon_{i_\ell} e_j)-u\|_{C^0(K)}=0
\end{equation}
for any compact $K\subset\overline\Omega_j\setminus \{p_j\}$ and any $j=1,...,N$. Also
\begin{equation}\label{spectral_stab_eig_3}
\lim_{\ell\to\infty}\|u_{\epsilon_{i_\ell}}-u\|_{H^1(\Omega_0)}=0\,,\ \ \lim_{\ell\to\infty}\|u_{\epsilon_{i_\ell}}-u\|_{C^0(K)}=0
\end{equation}
for all compact $K\subset\overline\Omega_0\setminus\{p_1,...,p_N\}$.
\end{theorem}
We remark that this theorem can be also deduced as an instance of the stability result \cite[Theorem 4.1 and Remark 4.3]{bucur_giacomini_trebeschi}.

The next result is a spectral stability result for the Steklov problem when a Lipschitz domain $\Omega_0$ is approximated by a family of Lipschitz (or smooth) domains $\{\Omega_\epsilon\}_{\epsilon\in[0,\epsilon_0)}$, as $\epsilon\to 0$ in a fixed atlas class, which is proved  in \cite[Theorem 2.5, Theorem 3.1 and Proposition 4.4]{FerreroLamberti}. We refer to \cite[Definition 2.6]{FerreroLamberti} for the precise definition of atlas class. Roughly speaking, a common atlas for a family $\{\Omega_\epsilon\}_{\epsilon\in[0,\epsilon_0)}$  is a fixed set of coordinate cylinders $\{C_i\}_{i=1}^m$, $C_i$ isometric to $\omega_i\times(0,L_i)$, in which each $\Omega_\epsilon$ is described, after a change of variables, as the subgraphs of a function $g_{i,\epsilon}:\omega_i\to\mathbb R$. The additional statement on the $C^0$ convergence up to the boundary follows from standard elliptic estimates. More precisely, it follows from boundary H\"older estimates for Robin problems on Lipschitz domains with uniformly bounded parameters, see \cite{nit}.

\begin{theorem}\label{thm:spec:stab}
Let $\Omega_0$, $\{\Omega_\epsilon\}_{\epsilon\in(0,\epsilon_0)}$ be bounded, Lipschitz domains in the same atlas class. Let $g_{i,\epsilon}:\omega_i\to\mathbb R$, $i=1,...,m$ denote the functions describing the boundary of $\Omega_\epsilon$ in each coordinate cylinder of the atlas. Assume that, for all $i=1,...,m$
$$
\lim_{\epsilon\to 0}\|g_{i,\epsilon}-g_{i,0}\|_{L^\infty(\omega_i)}= 0\,,\ \ \ \sup_{i=1,..,m\\\epsilon\in[0,\epsilon_0)}\|\nabla g_{i,\epsilon}\|_{L^{\infty}(\omega_i)}<\infty
$$
and that
$$
\lim_{\epsilon\to 0}|\partial\Omega_\epsilon|=|\partial\Omega_0|.
$$
Then, for all $k\in\mathbb N$
$$
\lim_{\epsilon\to 0}\sigma_k(\Omega_\epsilon)=\sigma_k(\Omega_0).
$$
Moreover, let $\{u_{\epsilon}\}_{\epsilon\in(0,\epsilon_0)}$ be a family of eigenfunctions associated with $\sigma_k(\Omega_\epsilon)$ normalised by $\int_{\partial\Omega_\epsilon}u_\epsilon^2=1$. Then, for any sequence $\epsilon_i\to 0$ as $i\to\infty$ there exists a subsequence $\epsilon_{i_\ell}\to 0$ as $\ell\to\infty$ and a Steklov eigenfunction $u$ associated with $\sigma_k(\Omega_0)$ normalised by $\int_{\partial\Omega_0}u^2=1$ such that
$$
\lim_{\ell\to\infty}\|u_{\epsilon_{i_\ell}}-u\|_{H^1(\Omega_0\cap\Omega_{\epsilon_{i_\ell}})}=0
$$
and along the same subsequence
$$
\lim_{\ell\to\infty}\|u_{\epsilon_{i_\ell}}-u\|_{C^0(\overline{\Omega_0\cap\Omega_{\epsilon_{i_\ell}}})}=0
$$
If $\sigma_k(\Omega_0)$ is simple, then $\sigma_k(\Omega_\epsilon)$ is simple for all $\epsilon>0$ small enough and the sign of $u_\epsilon$ can be chosen so  that convergence is along the whole family $\{\Omega_\epsilon\}_{\epsilon\in(0,\epsilon_0)}$.

\end{theorem}

Next, we recall a spectral stability result for the Steklov problem under tubular excisions in dimension $n\geq 3$. The original approach was developed in \cite{FS_SO} for the first eigenvalue, and the general result is proved in \cite[Theorem 1.1]{hong}.

\begin{theorem}\label{thm:excisions}
Let $\Omega\subset\mathbb R^n$, $n\geq 3$, be a bounded smooth domain. Let $\Gamma=\bigcup_{i=1}^m\gamma_i$ be the union of smooth, embedded, pairwise disjoint curves $\gamma_i:[0,1]\to\mathbb R^n$, $\gamma_i((0,1))\subset\Omega$ and meeting $\partial\Omega$ orthogonally at their extrema $\gamma_i(0),\gamma_i(1)$. Let $\Omega_\rho:=\{x\in\Omega :{\rm dist}(x,\Gamma)>\rho\}$. Then, for all $k\in\mathbb N$
$$
\lim_{\rho\to 0}\sigma_k(\Omega_\rho)=\sigma_k(\Omega).
$$
Moreover, if $\{u_\rho\}_{\rho\in(0,\rho_0)}$ is a family of eigenfunctions associated to $\sigma_k(\Omega_\rho)$ satisfying $\int_{\partial\Omega_\rho} u_\rho^2=1$, then for any sequence $\rho_i\to 0$ as $i\to\infty$ there exists a subsequence $\rho_{i_\ell}\to 0$ as $\ell\to\infty$ and a Steklov eigenfunction $u$ associated with $\sigma_k(\Omega)$ normalised by $\int_{\partial\Omega}u^2=1$ such that
$$
\lim_{\ell\to\infty}\|u_{\rho_{i_\ell}}- u\|_{C^0(K)}=0
$$
for all compact $K\subset\overline\Omega\setminus\Gamma$. If $\sigma_k(\Omega)$ is simple, the sign of $u_\rho$ can be chosen so that the convergence is along the whole family $\{\Omega_\rho\}_{\rho\in(0,\rho_0)}$.
\end{theorem}

We conclude this section by recalling the following generic simplicity result for Steklov eigenvalues \cite[Theorem 1.1]{wang_generic_euclid} (see also \cite{BPP_generic,wang_generic_manifold}).

\begin{theorem}\label{thm:genericity}
Let $\Omega$ be a smooth bounded domain in $\mathbb R^n$. Then for all $N\geq 1$, $m\geq 2$ and $\epsilon>0$ there exists a smooth embedding $\Phi:\overline\Omega\to\mathbb R^n$ satisfying
$$
\|\Phi-Id\|_{C^m(\overline\Omega)}<\epsilon
$$
such that $\sigma_1(\Phi(\Omega)),...,\sigma_N(\Phi(\Omega))$ are simple. The map $\Phi$ can be chosen to be the identity on $\overline \Omega\cap U$, for any prescribed open set $U\subset\mathbb R^n$ intersecting $\partial\Omega$ and such that $\partial\Omega\setminus\overline U\ne\emptyset$.
\end{theorem}
We remark that the last statement of Theorem \ref{thm:genericity} on the localisation of the perturbation is not explicitly stated in \cite[Theorem 1.1]{wang_generic_euclid}, however it is contained in the proof, in particular, in the proof of Proposition 4.14 of \cite{wang_generic_euclid}.

\section*{Acknowledgements}
The first author acknowledges support from the European Research Council (ERC) under the European Union's Horizon 2020 research and innovation programme through grant agreement~862342, from the Spanish State Research Agency through grants PID2022-136795NB-I00 and PID2025-168410NB-I00, and from the ICMAT--Severo Ochoa grant CEX2023-001347-S. The second author acknowledges support of the INdAM-GNAMPA group. The third author acknowledges support of the INdAM-GNSAGA group.

\section*{Use of AI tools}
ChatGPT-6.0 Astra was lightly used in the development and preparation of this work, including assistance with mathematical arguments and the revision of the manuscript.

\bibliography{bibliography}{}
\bibliographystyle{abbrv}

\end{document}